\documentclass[12pt]{article}

\usepackage[cp1251]{inputenc}  
\usepackage[T2A, OT1]{fontenc}
\usepackage[english]{babel}
\usepackage{amsfonts, longtable}
\usepackage{amsmath, amssymb, amsthm}
\usepackage{cite}
\usepackage{indentfirst}

\newtheorem{theo}{Theorem}[section]
\newtheorem{predl}[theo]{Statement}
\newtheorem{lem}[theo]{Lemma}
\newtheorem{sld}[theo]{Corollary}

\newcommand{\comment}[1]{}

\newenvironment{prf}{\par\textbf{Proof.}}{\par\qed\par}

\def\Ls{{{\mathcal L}}}

\title{Some Closed Classes of Three-Valued Logic Generated by Symmetric Functions}
\author{A.\,V.\,Mikhailovich\footnote{National Research University Higher School of Economics}}
\date{March, 2015}
\begin{document}
\maketitle
\begin{abstract}
Closed classes of three-valued logic generated by symmetric funtions that equal $1$ in almost all tuples from $\{1,2\}^n$ and equal $0$ on the rest tuples are considered. Criteria for bases existence for these classes is obtained.
\end{abstract}
The problem of the bases existence for some families of closed 
classes of the three-valued logic functions
is considered in the paper. 
E.\,Post~\cite{post41} (see also, for instance,~\cite{lau06}) described all closed classes of Boolean functions and showed that each
such class has a finite basis. This result is not extendable to the case of $k$-valued logics for 
$k\geq 3.$ Ju.\,I.\,Janov and A.\,A.\,Muchnik~\cite{janov} (see also, for instance,~\cite{lau06}) 
showed that for all $k\geq 3$ the set $P_k$ (here $P_k$ is the set of all functions of the $k$-valued logic) contains closed classes having a countable basis, and those having no basis. The generating
systems for classes from these examples consist of symmetric functions 
that take values from the set $\{0, 1\}$
and equal to zero on the unit tuples and on tuples containing at least one zero component.

Let $R$ be the set of all functions that take values from the set $\{0,1\}$
and equal zero on tuples containing at least one zero\footnote{All necessary definitions can be found in~\cite{lau06,mikh08,mikh09,mikh12}.}. 
\begin{multline*}
R=\{f(x_1,\ldots,x_n)\mid ((\forall\widetilde\alpha)((\widetilde\alpha\in\{0,1,2\}^n)
\rightarrow (f(\widetilde\alpha)\in\{0,1\}))) 
\& \\((\forall\widetilde\alpha)((\widetilde\alpha\in\{0,1,2\}^n\backslash\{1,2\}^n)
\rightarrow (f(\widetilde\alpha)=0)))\}
\end{multline*}
In this paper we deal with some subclasses of the class $R$.
Any function that does not change with variable relabeling is called symmetric.
Denote by $S$ the set of all symmetric functions from $R$.
The set of all tuples that can be obtained from each other by 
component permutation is called a layer. 
With $\mathcal L(e,d)$ we denote a layer from $\{1,2\}^n$
containing $e$ $1$s and $d$ $2$s.
Let $N_f=\{\widetilde\alpha\in E_3^n\mid f(\widetilde\alpha)=1\}.$
With $i_s(x_1,\ldots,x_s)$ we denote a function from $R$ 
such that $N_{i_s}=\{1,2\}^s.$
Let 
$$
S_{-1}=\{f(x_1,\ldots,x_n)\in R\mid (\exists e) (\exists d) 
(N_f=\{1,2\}^n\backslash\mathcal L(e,d))\}.
$$
Let $f\in S_{-1}.$
Let us denote by $e_f$ and $d_f$ such numbers that
$N_f=\{1,2\}^n\backslash\mathcal L(e_f,d_f).$

\begin{predl}
Let $\Phi$ be a formula over $R;$
$\Phi_1$ be a subformula of $\Phi.$
Then for any tuple $\widetilde\alpha$ 
the equality $\Phi_1(\widetilde\alpha)=0$ implies
the equality $\Phi(\widetilde\alpha)=0.$
\label{mikh_mos_lemma10}
\end{predl}
\begin{prf}
Suppose $d$ is the depth of $\Phi,$ $d_1$ is the depth of $\Phi_1.$
Let $\Phi$ have the form 
$g(\mathcal B_1,\ldots,\mathcal B_m).$ 

The proof is by induction on $d-d_1.$
For $d-d_1=0,$ the proof is trivial.
Let $d-d_1=1.$ 
Then $\Phi_1=\mathcal B_t,$ $1\leq t\leq m.$
Since $g\in R,$ 
$\mathcal B_t(\widetilde\alpha)=0$ implies $\Phi(\widetilde\alpha)=0.$

Let $d-d_1>1.$
Then $\Phi_1$ is a subformula of $\mathcal B_t$ for $1\leq t\leq m.$
Denote by $d_2$ the depth of $\mathcal B_t.$
Then $d_2-d_1<d-d_1.$
By the inductive assumtion $\mathcal B_t(\widetilde\alpha)=0.$
Note that $d-d_2<d-d_1.$
Then by the inductive assumtion $\Phi(\widetilde\alpha)=0.$
\end{prf}

\begin{sld}
If a formula over $R$ equals $1$ on a tuple$,$ 
then any subformula of the formula also equals $1$ on the same tuple.
\label{mikh_mos_sld1010}
\end{sld}

\begin{sld}
If $\Phi$ is a formula over $R$ and 
$\Phi_1$ is a subformula of formula $\Phi$ then $N_{\Phi}\subset N_{\Phi_1}.$
\label{mikh_mos_sld1020}
\end{sld}

\begin{predl}
If $f(x_1,\ldots,x_n)\in S_{-1},$ $n>3,$ $e_f>0,$ $d_f>0,$ 
then $i_2(x_1,x_2)\in[\{f\}].$
\label{mikh_mos_lemma20}
\end{predl}
\begin{prf}
Suppose $e_f>d_f>1.$ 
Then 
$i_2(x_1,x_2)=f(\underbrace{x_1,\ldots,x_1}_{e_f+1},\underbrace{x_2,\ldots,x_2}_{d_f-1}).$

Suppose $d_f>e_f>1.$ 
Then 
$i_2(x_1,x_2)=f(\underbrace{x_1,\ldots,x_1}_{d_f+1},\underbrace{x_2,\ldots,x_2}_{e_f-1}).$

Suppose $e_f=1.$ Then $d_f>2$ and
$i_2(x_1,x_2)=f(x_1,x_1,\underbrace{x_2,\ldots,x_2}_{d_f-1}).$

Suppose $d_f=1.$ Then $e_f>2$ and
$i_2(x_1,x_2)=f(x_1,x_1,\underbrace{x_2,\ldots,x_2}_{e_f-1}).$
\end{prf}

\begin{sld}
If $f(x_1,\ldots,x_n)\in S_{-1},$ $n>3,$ $e_f>0,$ $d_f>0$ then $i_s(x_1,\ldots,x_s)\in[\{f\}]$ 
for any $s\in{\mathbb N}.$
\label{mikh_mos_sld2010}
\end{sld}

\begin{predl}
Let $f(x_1,\ldots,x_n)\in S_{-1};$ 
$f$ be obtained by means of formula $\Phi$ over $S_{-1};$
$g(\mathcal B_1,\ldots,\mathcal B_m)$ be a subformula of $\Phi;$  $m<n.$ 
Then $g(\mathcal B_1,\ldots,\mathcal B_m)=i_m(\mathcal B_1,\ldots,\mathcal B_m).$
\label{mikh_mos_lemma30}
\end{predl}
\begin{prf}
Without loss of generality $\mathcal B_1,\ldots,\mathcal B_s$ are non-trivial
formulae and $\mathcal B_{s+1},\ldots,\mathcal B_m$ are variables 
$x_1,$\ldots,$x_t.$
Obviously $t\leq m<n.$
Let us denote $g(\mathcal B_1,\ldots,\mathcal B_m)$ by $\Phi_1,$
$i_m(\mathcal B_1,\ldots,\mathcal B_m)$ by $\Phi_2.$

Let $\widetilde\alpha=(\alpha_1,\ldots,\alpha_n)\in N_f;$
$\beta_i=\mathcal B_i(\widetilde\alpha)$ $(i=1,\ldots,m),$
$\widetilde\beta=(\beta_1,\ldots,\beta_m).$ 
Suppose $\widetilde\beta\in\Ls(e_g, d_g).$
Consider tuples $\widetilde\alpha^0,$ \ldots, $\widetilde\alpha^{n-t}$ such that
$\widetilde\alpha^j=(\alpha_1,\ldots,\alpha^t,1^j,2^{n-t-j}).$
Since $n>t$ there exists $j,$ $0\leq j\leq n-t$,
such that $\widetilde\alpha^j\in N_f.$ 
Using Corollary~\ref{mikh_mos_sld1010} we get 
$\mathcal B_i(\widetilde\alpha^j)=1$ $(i=1,\ldots,s),$ $\Phi_1(\widetilde\alpha^j)=1.$
Since $\alpha_i=\alpha_i^j,$ $(i=1,\ldots,t),$
$\mathcal B_i(\widetilde\alpha^j)=\mathcal B_i(\widetilde\alpha).$ 
Thus, $g(\widetilde\beta)=1.$
That contradicts $\widetilde\beta\in\Ls(e_g,d_g).$
Therefore, $\widetilde\beta\in N_g.$
Hence,
$\Phi_1(\widetilde\alpha)=\Phi_2(\widetilde\alpha).$
\end{prf}

\begin{sld}
Let $f(x_1,\ldots,x_n)\in S_{-1};$ 
$f$ be obtained by means of formula $\Phi$ over $S_{-1}.$
Then there is subformula $g(\mathcal B_1,\ldots,\mathcal B_m)$ of the formula $\Phi$
such that $m\geq n.$
\label{mikh_mos_sld3010}
\end{sld}

\begin{predl}
Suppose $\mathfrak A\subseteq S_{-1},$ $f(x_1,\ldots,x_n)\in S_{-1},$ 
$n>3,$ $e_f>0,$
$f\in[\mathfrak A]$.
Then there is function $g\in\mathfrak A$ such that $f\in[\{g\}].$
\label{mikh_mos_lemma40}
\end{predl}
\begin{prf}
Let $f(x_1,\ldots,x_n)$ be obtained by means of formula $\Phi$ over $\mathfrak A;$
$\Psi$ be a subformula of $\Phi$ of the form
$g(\mathcal B_1,\ldots,\mathcal B_m);$
$N_{\Psi}\neq \{1,2\}^n$ and any 
non-trivial subformula $h(\mathcal C_1,\ldots,\mathcal C_q)$
of the formula $\Psi$ equal
$i_q(\mathcal C_1,\ldots,\mathcal C_q).$
Hence, $g(\mathcal B_1,\ldots,\mathcal B_m)$
cannot be replaced by $i_m(\mathcal B_1,\ldots,\mathcal B_m).$
Thus, $m\geq n>3.$
Let $p=n!;$ 
$\widetilde x=(x_1,\ldots,x_n),$ $\sigma$ be a permutation of $\{1,\ldots,n\};$
Denote $\sigma(\widetilde x)=(x_{\sigma(1)},\ldots,x_{\sigma_n}).$
Let $\sigma_1,$\ldots,$\sigma_p$ be a sequence of all permutations of $\{1,\ldots,n\}.$

Consider 2 cases. 

Case 1: $\Psi$ is a simple formula.
Consider 3 subcases. 

Subcase 1: $N_\Psi=N_f.$ Then obviously $f\in[\{g\}].$

Subcase 2: $N_\Psi\subset  N_f,$ $N_\Psi\neq N_f,$ $e_g>0,$ $d_g>0.$
Using Corollary~\ref{mikh_mos_sld2010} we get
$i_p\in[\{g\}].$ 
Then
$$
f(x_1,\ldots,x_n)=i_p(\Psi(\sigma_1(\widetilde x)),\ldots,\Psi(\sigma_p(\widetilde x))).
$$

Subcase 3: $N_\Psi\subset  N_f,$ $e_f=0$ or $d_f=0.$ Then 
$N_\Psi=N_f$ and $f\in[\{g\}].$

Case 2: $\Psi$ is not a simple formula.
Then $m>n>3$ and $e_g>0.$ 
Consider 2 subcases. Subcase 1: $d_g>0.$
Then $i_1, i_p\in [\{g\}].$
Let 
$$
\mathcal C_j=\left\{\begin{array}{l}
\mathcal B_j \mbox{ if } \mathcal B_j\in\{x_1,\ldots,x_n\}, \\
i_1(x_1) \mbox{ otherwise,}
\end{array}\right.
$$

$$
\Psi_1=g(\mathcal C_1,\ldots,\mathcal C_m).
$$

Obviously $\Psi_1$ is a formula over $\{g\}.$
If $N_{\Psi}=N_f$ then
$f(x_1,\ldots,x_n)=\Psi_1$ and $f\in[\{g\}].$
If $N_{\Psi}\neq N_f$ then
$$
f(x_1,\ldots,x_n)=i_p(\Psi_1(\sigma_1(\widetilde x)),\ldots,\Psi_1(\sigma_p(\widetilde x))).
$$

Subcase 2: $d_g=0.$ Then $d_f=0.$ Hence
$$
f(x_1,\ldots,x_n)=g(\underbrace{x_1,\ldots,x_1}_{m-n+1},x_2,x_3,\ldots,x_n).
$$
Thus, $f\in[\{g\}].$
\end{prf}

\begin{predl}
Suppose $\mathfrak A\subseteq S_{-1},$ $f(x_1,\ldots,x_n)\in S_{-1},$ 
$n>3,$ $e_f=0,$
$f\in[\mathfrak A]$.
Then there is a function $g\in\mathfrak A$ such that $f\in[\{g\}].$
\label{mikh_mos_lemma50}
\end{predl}
\begin{prf}
Let $f(x_1,\ldots,x_n)$ be obtained by means of formula $\Phi$ over $\mathfrak A.$
Consider 2 cases. 

Case 1: there is a subformula of $\Phi$ of the form $g(x_{i_1},\ldots,x_{i_m})$ such that $m\geq n.$ 
Then $e_g=0$ and
$$
f(x_1,\ldots,x_n)=g(\underbrace{x_1,\ldots,x_1}_{m-n+1},x_2,x_3,\ldots,x_n).
$$
Thus, $f\in[\{g\}].$

Case 2: there is subformula of $\Phi$ of the form $g(x_{i_1},\ldots,x_{i_m})$ such that $m<n.$
Without loss of generality let $\{x_{i_1},\ldots,x_{i_m}\}=\{x_1,\ldots,x_s\},$ $s\leq m.$
Consider tuple $\widetilde\alpha=(1^s,2^{n-s})\in N_f.$ 
Obviously $g(\alpha_{i_1},\ldots,\alpha_{i_m})=0.$ Thus, $f(x_1,\ldots,x_n)$ 
cannot be obtained by means of formula $\Phi.$ Hence, case 2 is impossible.
\end{prf}

\begin{lem}
Let $F\subset S_{-1};$ $F$ not contain
congruent functions$;$ $F$ be infinite set$;$ for any $e\in\mathbb N,$
$d\in \mathbb N$ there be function $h\in F$ such that 
$e_h>e,$ $d_h>d;$
and if $f\in F$ then $e_f>0$ and $d_f>0$.  
Then for any $f\in F$ there is function 
$g\in F$ such that $f\in[\{g\}].$
\label{mikh_mos_lemma54}
\end{lem}
\begin{prf}
Let $f(x_1,\ldots,x_n)\in F;$
$g\in F$ such that $e_g\geq 2^{e_f}-1,$
$d_g\geq 2^{e_f}(2^{d_f}-1)$. Using Corollary~\ref{mikh_mos_sld2010} we obtain
that $i_1\in [\{g\}],$ $i_p\in[\{g\}].$ Let
$$
f_1=g(\underbrace{i_1(x_1),\ldots,i_1(x_1)}_{e_g-2^{e_f}+1},
x_1,\ldots,\underbrace{x_j,\ldots,x_j}_{2^{j-1},}\ldots,
\underbrace{x_{n-1},\ldots,x_{n-1}}_{2^{n-2}},
\underbrace{x_n,\ldots,x_n}_{d_g-2^n+2^{e_f}}). 
$$
Let $p=n!$.
Let $\widetilde x=(x_1,\ldots,x_n),$ $\sigma$ be permutation of $\{1,\ldots,n\}.$
Denote $\sigma(\widetilde x)=$ $(x_{\sigma(1)},\ldots,x_{\sigma_n}).$
Let $\sigma_1,$\ldots,$\sigma_p$ be a sequence of all permutations of $\{1,\ldots,n\}.$
It is clear that $(1^{e_f},2^{d_f})\notin N_{f_1}$ and 
if $\widetilde\alpha\in\{1,2\}^n,$ and $\widetilde\alpha\neq (1^{e_f},2^{d_f})$ 
then $\widetilde\alpha\in N_{f_1}.$ Hence,
$$
f(x_1,\ldots,x_n)=
i_p(f_1(\sigma_1(\widetilde x)),\ldots,f_1(\sigma_p(\widetilde x))).
$$
Thus, $f\in[\{g\}].$
\end{prf}

\begin{sld}
Let $i=1,2;$ $F_i\subset S_{-1};$ $F_i$ not contain
congruent functions$;$ $F_i$ be infinite set$;$ for any $e_i\in\mathbb N,$
$d_i\in \mathbb N$ there exist function $h\in F$ such that 
$e_h>e_i,$ $d_h>d_i;$
and if $f\in F_i$ then $e_f>0$ and $d_f>0$.  
Then $[F_1]=[F_2].$
\label{mikh_mos_sld5410}
\end{sld}

\begin{sld}
There are countably many classes of the form $G=[F]$ such
that $F\subset S_{-1}$ and for any $e,d\in\mathbb N$ there exists a function
$h\in F$ such that $e_h>e,$ $d_h>d.$
\label{mikh_mos_sld5420}
\end{sld}

\begin{lem}
Let $F\subset S_{-1};$ $F$ not contain
congruent functions$;$ $F$ be infinite set$;$ there exist $e\in\mathbb Z^{+}$
such that if $f(x_1,\ldots,x_n)\in F$ then $e_f=e.$
Then for any $f\in F$ there exists a function $g\in F$ such that $f\in [\{g\}].$
\label{mikh_mos_lemma55}
\end{lem}
\begin{prf}
Let $f(x_1,\ldots,x_n)\in F.$ 
If $e=0$, $m>n$ and $g(x_1,\ldots,x_m)\in F$ then
$$
f(x_1,\ldots,x_n)=g(\underbrace{x_1,\ldots,x_1}_{m-n+1},x_2,x_3,\ldots,x_n).
$$
Consider $e>0.$
Suppose $g(x_1,\ldots,x_m)\in F$ such that $d_g>(e_f+1)d_f.$
Using Corollary~\ref{mikh_mos_sld2010} we obtain
that $i_1\in [\{g\}],$ $i_p\in[\{g\}].$ Let
$$
f_1=g(x_1,\ldots,x_{e_f}\underbrace{x_{e_f+1},\ldots,x_{e_f+1}}_{e_f+1},
\ldots,\underbrace{x_{n-1},\ldots,x_{n-1}}_{e_f+1},
\underbrace{x_n,\ldots,x_n}_{d_g-(e_f+1)(d_f-1)}).
$$
Let $p=n!$.
Let $\widetilde x=(x_1,\ldots,x_n),$ $\sigma$ be permutation of $\{1,\ldots,n\}.$
Denote $\sigma(\widetilde x)=$ $(x_{\sigma(1)},\ldots,x_{\sigma_n}).$
Let $\sigma_1,$\ldots,$\sigma_p$ be a sequence of all permutations of $\{1,\ldots,n\}.$
It is clear that $(1^{e_f},2^{d_f})\notin N_{f_1}.$ Furthermore, 
if $\widetilde\alpha\in\{1,2\}^n$ and $\widetilde\alpha\neq (1^{e_f},2^{d_f})$ 
then $\widetilde\alpha\in N_{f_1}.$ Hence,
$$
f(x_1,\ldots,x_n)=
i_p(f_1(\sigma_1(\widetilde x)),\ldots,f_1(\sigma_p(\widetilde x))).
$$
Thus, $f\in[\{g\}].$
\end{prf}

\begin{sld}
Let $F_1, F_2\subset S_{-1};$ $F_1$ and $F_2$ not contain
congruent functions$;$ $F_1,$ $F_2$ be infinite$;$ there exist $e\in\mathbb Z^+$
such that if $f(x_1,\ldots,x_n)\in F_1\cup F_2$ then $e_f=e$.
Then $[F_1]=[F_2].$
\label{mikh_mos_sld5510}
\end{sld}

\begin{sld}
There are countably many classes of the form $G=[F]$ such
that there exists $e\in\mathbb Z^{+}$
such that if $f(x_1,\ldots,x_n)\in F$ then $e_f=e$.
\label{mikh_mos_sld5520}
\end{sld}

\begin{lem}
Let $F\subset S_{-1};$ $F$ not contain
congruent functions, $F$ be infinite$;$ there exist $d$
such that if $f(x_1,\ldots,x_n)\in F$ then $d_f=d$.
Then for any $f\in F$ there exists a function $g\in F$ such that $f\in [\{g\}].$
\label{mikh_mos_lemma56}
\end{lem}
The proof is similar to the proof of Lemma~\ref{mikh_mos_lemma55}.
\begin{sld}
Let $F_1, F_2\subset S_{-1};$ $F_1$ and $F_2$ not contain
congruent functions$;$ $F_1,$ $F_2$ be infinite$;$ there exist $d$
such that if $f(x_1,\ldots,x_n)\in F_1\cup F_2$ then $d_f=d$.
Then $[F_1]=[F_2].$
\label{mikh_mos_sld5610}
\end{sld}

\begin{sld}
There are countably many classes of the form $G=[F]$ such
that there exist $d\in\mathbb Z^{+}$
such that if $f(x_1,\ldots,x_n)\in F$ then $d_f=d$.
\label{mikh_mos_sld5620}
\end{sld}

\begin{lem}
Let $F\subset S_{-1};$ $F$ not contain
congruent functions$;$ $F$ be infinite. Then there exists $N$
such that for any function $f(x_1,\ldots,x_n)$ with $n>N$ there exists
a function $g\in F$ such that $f\in[\{g\}].$
\label{mikh_mos_lemma60}
\end{lem}
\begin{prf}
Let $F_0=\{f\in F\mid e_f>0, d_f>0\},$
$F_1=\{f\in F\mid e_f=0\},$ $F_2=\{f\in F\mid d_f=0\}.$
Obviously $F=F_0\cup F_1\cup F_2.$
Let
$$
N_j=\left\{\begin{array}{l}
\max_{f(x_1,\ldots,x_n)\in F_j} (n), \mbox{ if } |F_j|<\infty; \\
0, \mbox{ otherwise,}
\end{array}
\right.
$$
$j=1,2.$

Let $F_0$ be infinite and $f(x_1,\ldots,x_n)\in F_0.$ 
Consider 3 cases. Case 1: for any $e\in \mathbb N,$ $d\in \mathbb N,$
there exists $h\in F_0$ such that  $e_h>e$ and $d_h>d$. 
Let $N_0=1,$ $N=\max(N_0, N_1, N_2).$
Using Lemma~\ref{mikh_mos_lemma54} we obtain that there is
function $g\in F_{e_f}$ such that $f\in[\{g\}].$

Case 2: there exists $e_0\in \mathbb N$ such that if
$h\in F_0$ then $e_h\leq e_0$.
Let
$F^s=\{f\in F_0\mid e_f=s\};$
$$
n_s=\left\{\begin{array}{l}
\max(n)_{f(x_1,\ldots,x_n)\in F^s } \mbox{ if } |F^s|<\infty; \\
0 \mbox{ otherwise, } 
\end{array}\right.
$$
$s=1,\ldots,e_0.$
Let  $N_0=\max(n_1,\ldots,n_s),$ $N=\max(N_0, N_1, N_2).$
Since $e_f>n$ then $|F_{e_f}|=\infty.$ 
Using Lemma~\ref{mikh_mos_lemma55} we obtain that there is
a function $g\in F_{e_f}$ such that $f\in[\{g\}].$

Case 3: there exists $d_0\in \mathbb N$ such that if
$h\in F_0$ then $d_h\leq d_0$. This case is similar to case 2
with reference to Lemma~\ref{mikh_mos_lemma56} respectively.

Let $F_1$ be infinite and $f(x_1,\ldots,x_n)\in F_1.$ 
Then for any $g(x_1,\ldots,x_m)\in F_1$ such that $m>n$ we have
$$
f(x_1,\ldots,x_n)=g(\underbrace{x_1,\ldots,x_1}_{m-n+1},x_2,\ldots,x_n).
$$ 
The case of infinite $F_2$ is considered in the same way.

Thus for any $f(x_1,\ldots,x_n)\in F$ 
such that $n\geq N$ there exist $g\in F\backslash\{f\}$ such that 
$f\in[\{g\}].$
\end{prf}

\begin{lem}
$S\subset [S_{-1}].$
\label{mikh_mos_lemma70}
\end{lem}
\begin{prf}
Let $f\in S.$ Denote
$$
\overline N_f=\{\widetilde\alpha\in\{1,2\}^n\mid f(\widetilde\alpha)=0\}.
$$
Let $f_1,\ldots,f_t$ be the set of all functions from $S_{-1}$ 
such that $\Ls(e_{f_j}, d_{f_j})\in \overline N_f,$ $j=1,\ldots,t.$
Using Corollary~\ref{mikh_mos_sld2010} we obtain that $i_t\in S_{-1}.$
It is clear that
$$
f(x_1,\ldots,x_n)=i_t(f_1(x_1,\ldots,x_n),\ldots,f_t(x_1,\ldots,x_n)).
$$
Hence, $f\in[S_{-1}].$ Thus, $S\subset [S_{-1}].$
\end{prf}

\begin{lem}
$R\subset[S].$
\label{mikh_mos_lemma80}
\end{lem}
\begin{prf}
Let $f(x_1,\ldots,x_n)\in R;$ $N_f=\{\widetilde\alpha^1,\ldots,\widetilde\alpha^t\};$
$m=2^n-1;$
$$
e_i=\sum_{\alpha^i_j=1}2^{j-1}, d_i=m-e_i, i=1,\ldots,t.
$$
Let $g(x_1,\ldots,x_m)\in S$ such that
$N_g=\cup_{i=1}^t \Ls(e_i,d_i).$
It is easily shown that 
$$
f(x_1,\ldots,x_n)=g(x_1,x_2,x_2,\underbrace{x_3,\ldots,x_3}_4,\ldots,\underbrace{x_i,\ldots,x_i}_{2^{i-1}},
\ldots,\underbrace{x_n,\ldots,x_n}_{2^{n-1}}).
$$
Thus, $f\in[S].$ Hence, $R\subseteq [S].$
\end{prf}

\begin{theo}
$[S_{-1}]=R.$
\end{theo}
The statement of the theorem is a consequence of Lemma~\ref{mikh_mos_lemma70} and Lemma~\ref{mikh_mos_lemma80}.


\begin{theo}
Let $F\subset S_{-1};$ $F$ not contain
congruent functions$;$ $G=[F].$ Then
\begin{enumerate}
\item Class $G$ has a finite basis iff $F$ is finite.
\item Class $G$ has no basis iff $F$ is infinite.
\end{enumerate}

\end{theo}
\begin{prf}
1. If $F$ is finite, the proof is trivial. Let $G$ have a finite basis 
$\mathfrak A=\{f_1(x_1,\ldots,x_{n_1}),$ \ldots, $f_s(x_1,\ldots,x_{n_s})\}.$
Suppose 
$n_a=\max_{1\leq j\leq s} n_j.$ 
Since $G=[F]$ then $F\in[\mathfrak A].$ 
Using Corollary~\ref{mikh_mos_sld3010} we obtain that if $f(x_1,\ldots,x_n)\in F$ then
$n\leq n_a.$ Hence, $F$ is finite.

2. If $G$ has no basis, obviously $F$ is infinite. Let $F$ be infinite.
By $N$ denote such number that for any $f(x_1,\ldots,x_n)\in F$ 
such that $n>N$ there exist a function $g\in F$ such that $f\in[\{g\}].$
According to Lemma~\ref{mikh_mos_lemma60} there exist such a number.
Suppose that $G$ has basis $\mathfrak A.$ 
Let $F_N=\{f(x_1,\ldots,x_n)\in F\mid n\leq N\},$
$\mathfrak A_N\subset \mathfrak A$ such that $F_N\subset[\mathfrak A_N]$
and for any $\mathfrak A'\subset \mathfrak A_N,$ $\mathfrak A'\neq \mathfrak A_N$
$F_N\not\subset[\mathfrak A'].$
It is clear that $F_N$ is finite. Hence, $\mathfrak A_N$ is finite.

Consider $f\in\mathfrak A\backslash\mathfrak A_N.$
Let $\Phi$ be a formula over $F$ that realizes $f;$ and
$g_1,\ldots,g_s$ be all the functions from $F$ that are used in $\Phi.$
Let us prove that $g_j\in[\mathfrak A\backslash\{f\}],$ $j=1,\ldots,s.$
If $g_j\in F_N$ then $g_j\in[\mathfrak A_N].$ Hence, 
$g_j\in [\mathfrak A \backslash\{f\}].$
If $g_j\in F\backslash F_N$ there exist a function $h\in F\backslash (F_N\cup\{g_j\})$
such that $g_j\in[\{h\}].$ Hence, $f\in[\mathfrak A\backslash\{f\}].$
That contradicts the basis definition. Thus, $G$ has no basis.
\end{prf}

\begin{theo}
There are countably many classes of the form $G=[F]$ such 
that $F\subseteq S_{-1}.$
\end{theo}
The proof is based on Corollaries~\ref{mikh_mos_sld5420}, \ref{mikh_mos_sld5520}, 
and \ref{mikh_mos_sld5620}

\medskip 

This study (research grant No 14-01-0144) supported by The National Research University~--- Higher School of Economics' Academic Fund Program in 2014/2015.

\end{document}